\documentclass[11pt, english, reqno]{amsart}

\usepackage[T1]{fontenc}
\usepackage{babel}
\usepackage{mathrsfs}
\usepackage{amsthm}
\makeindex

\usepackage[dvips,top=3.8cm,left=4cm,right=3.8cm,
foot=3.8cm,bottom=4.0cm]{geometry}
\usepackage{amsfonts,amssymb,amsmath,amsthm, booktabs, 
latexsym}
\usepackage{centernot}
\usepackage{tikz-cd}
\usepackage{bbm}

\usepackage{enumerate}
\usepackage{color}

\newtheorem{theorem}{Theorem}[section]
\newtheorem{lemma}[theorem]{Lemma}
\newtheorem{proposition}[theorem]{Proposition}
\newtheorem{corollary}[theorem]{Corollary}
\newtheorem{question}[theorem]{Question}
\theoremstyle{definition}
\newtheorem{definition}[theorem]{Definition}

\theoremstyle{remark}
\newtheorem{remark}[theorem]{Remark}

\numberwithin{equation}{section}

\numberwithin{equation}{section}
\usepackage[pagewise]{lineno}

\newcommand{\be}{\begin{equation}}
\newcommand{\ee}{\end{equation}}
\newcommand{\ba}{\begin{aligned}}
\newcommand{\ea}{\end{aligned}}

\newcommand{\N}{{\mathbb N}}
\newcommand{\Z}{{\mathbb Z}}
\newcommand{\R}{{\mathbb R}}

\newcommand{\Leb}{\ensuremath{\lambda}}

\newcommand{\LS}{\ensuremath{\underset{n=1}{\overset{\infty}{\cap}} \, {\underset{i=n}{\overset{\infty}{\cup}}}\,}}

\def\csi1{\circ\sigma^{-1}}

\title[A shrinking target theorem for ergodic transformations of the unit interval] 
      {A shrinking target theorem for ergodic transformations of the unit interval}

\author[Shrey Sanadhya]{}

\subjclass[2020]{37A05, 37A10, 37E10}
 \keywords{Shrinking target problem, Interval exchange transformations}

 \email{shreysanadhya@gmail.com}

\begin{document}
\maketitle

\centerline{\scshape Shrey Sanadhya}
\medskip
{\footnotesize
 \centerline{Department of Mathematics}
   \centerline{Ben Gurion University of the Negev}
   \centerline{Be’er Sheva, 8410501, Israel.}
} 

\bigskip

\begin{abstract} We show that for \textit{any} ergodic Lebesgue measure preserving transformation $f: [0,1) \rightarrow [0,1)$ and \textit{any} decreasing sequence $\{b_i\}_{i=1}^{\infty}$ of positive real numbers with divergent sum, the set $$\LS f^{-i}(B (R_{\alpha}^{i} x,b_i))$$ has full Lebesgue measure for almost every $x \in [0,1)$ and almost every $\alpha \in [0,1)$. Here $B(x,r)$ is the ball of radius $r$ centered at $x \in [0,1)$ and $R_{\alpha}: [0,1) \rightarrow [0,1)$ is rotation by $\alpha \in [0,1)$. As a corollary, we provide partial answer to a question asked by Chaika (Question $3$, \cite{Chaika_2011}) in the context of interval exchange transformations. 

\end{abstract}

\section{ Introduction}\label{sect Intro} In this paper, we study a \textit{shrinking target problem} in the context of ergodic Lebesgue measure preserving transformations of the unit interval. In particular, we consider the case when the shrinking targets are centered along the orbit of an irrational rotation.

A typical shrinking target problem involves a probability measure preserving dynamical system $(X,\mu, T)$ and a sequence of measurable sets $\{B_n\}_{n \in \N}$ such that $\mu(B_n) \rightarrow 0$ (which are called \textit{shrinking targets}). We ask if  $$\mu (\{x \in X : T^n x \in B_n \,\, \textrm{for infinitely many}\,\, n\in \mathbb{N} \}) = 1$$ for $\mu$-a.e. $x \in X$. The Borel-Cantelli lemma gives us the necessary condition i.e. $\underset{n=1}{\overset{\infty}{\sum}} \mu(B_n)=\infty$. Hence in the setting of a metric space, shrinking balls with fixed centers and divergent radii sum, form a natural starting point for such a study. One of the earliest result in this setting is due to J. Kurzweil \cite{Kurzweil_1955}, who studied the shrinking target problem in the context of irrational rotations. 

For $\alpha \in [0,1)$, we define $R_{\alpha} \colon [0,1) \to[0,1)$ as $R_{\alpha}(x)=x+\alpha-\lfloor x+\alpha \rfloor$. Let $B(x,r)$ be the ball of radius $r$ centered at $x \in [0,1)$ and $\lambda$ be the Lebesgue measure on $[0,1)$. Kurzweil proved the following result :

\begin{theorem}[J. Kurzweil \cite{Kurzweil_1955}]  For any decreasing sequence of positive real numbers $\{b_i\}_{i=1}^{\infty}$ with divergent sum and for almost every $\alpha \in [0,1)$,

$$\Leb\left(\LS B(R_{\alpha}^i(x),b_i)\right)=1 $$ every $x \in [0,1)$.
\end{theorem} The pioneering work of Kurzweil motivated the study of shrinking target problem in diverse settings. We refer interested readers to \cite{Chernov_Kleinbock_2001},\cite{Maucourant_2006},\cite{Kim_2007},\cite{Galatolo_Kim_2007},\cite{Fernandez_Melian_Pestana_2012} and the survey paper \cite{Athreya_2009}. 

J. Chaika \cite{Chaika_2011} studied shrinking target problem in context of interval exchange transformations (IETs) (see Definition \ref{IET}). He showed that as a family, interval exchange transformations satisfy the shrinking target property (see Theorem \ref{Chaika}). This paper is motivated by results in \cite{Chaika_2011}. In particular, we provide a partial answer to Question $3$ of \cite{Chaika_2011}. We mention the question below for completion.

\begin{question}\label{q1} Let $\{y_i\}_{i=1}^{\infty} $ be a sequence of points in $ [0,1)$ and $\{b_i\}_{i=1}^{\infty}$ be a sequence of positive real numbers with divergent sum. Is it true that for almost every IET $T$, we have $\Leb(\LS T^{-i}(B(y_i,b_i)))=1$?
\end{question}

We answer the above question for the case when the sequence $\{y_i\}$ is an orbit under irrational rotations. To our surprise, if the shrinking targets are centered at points along the orbit of an irrational rotation, Question $\ref{q1}$ can be answered in much more generality, i.e. for \textit{any} ergodic Lebesgue measure preserving transformation of the unit interval and \textit{any} sequence $\{b_i\}_{i=1}^{\infty}$ of positive real numbers with divergent sum (see Theorem \ref{Main}). Thus, we obtain a partial answer to Chaika's question as a corollary (see Corollary \ref{Cor}). 

The \textit{outline} of the paper is as follows. In 
Section \ref{Prelim}, we provide background and state the main results (Theorem \ref{Main} and Corollary \ref{Cor}). Section \ref{proof} is dedicated to the proof of Theorem \ref{Main}.

\section{Preliminaries and statement of results}\label{Prelim}

\begin{definition}\label{IET} Let $P=(p_1,p_2,...,p_d)$, where $p_i \geq 0$ for each $i \in \{1,2,...,d\}$ be a $d$-dimensional vector which defines interval $ I = [0,\underset{i=1}{\overset{d}{\sum}} p_i)$, with $d$ sub intervals $$I_1 = [0,p_1), I_2 = [p_1, p_1+p_2),....,I_d = [p_1+...+p_{d-1}, p_1+...+p_{d-1}+p_d).$$  Let $\pi$ be a fixed permutation on the set $\{1,2,...,d\}$. An \textit{interval exchange transformation} (IET) is a map $T : I \rightarrow I$ which permutes the $d$ sub intervals $I_i$ by $\pi$. In other words for $x \in I_j$  $$T(x)= x - \underset{t<j}{\sum} p_t +\underset{\pi(t')<\pi(j)}{\sum} p_{t'}.$$ 
\end{definition} The IET defined above is on $d$-intervals and will be referred to as a $d$-IET. In this paper we will restrict ourselves to IETs on unit the interval, i.e. $I = [0,1)$ for convenience. A permutation $\pi \in S_d$ is called \textit{irreducible} if $\pi(\{1,...,t\}) \neq \{1,...,t\}$ for any $t < d$. It was shown by M. Keane \cite{Keane_1975} that IETs with dense orbits have irreducible permutations. Hence these IETs are important from the point of view of shrinking target property. In this paper we will work with $d$-IETs (for $d > 2$) with a fixed irreducible permutation $\pi \in S_d$ on the unit interval $I=[0,1)$. We can parametrize such IETs by the $d$-dimensional simplex $\Delta^d = \{(p_1,p_2,...,p_d) : p_i\geq 0, \underset{i}{\sum}\, p_i = 1\}$. 

Note that the standard simplex $\Delta^d$ comes equipped with Lebesgue measure. In this paper for a fixed irreducible permutation, the term almost every IET will refer to the Lebesgue measure on simplex $\Delta^d$. We will consider it as the parameterizing space of IETs in this paper. S. Kerkchoff, H. Masur, and J. Smillie (see \cite{Kerckhoff_Masur_Smillie_1986}) showed that for almost every IET is uniquely ergodic with respect to $\lambda$ (Lebesgue measure on the interval $[0,1)$). Below we list some known shrinking target results for IETs.

\begin{theorem} [J. Chaika \cite{Chaika_2011}]\label{Chaika} For almost every IET $T$ and a decreasing sequence $\{b_i\}_{i=1}^{\infty}$ of real numbers with divergent sum,

\begin{center}
$\LS T^{-i}\left(B(y,b_i)\right)$ has full $\lambda$-measure for every
$y$.

\end{center}

\end{theorem} Following logarithm law is due to S. Galatolo \cite{Galatolo_2006}.

\begin{theorem}[S. Galatolo \cite{Galatolo_2006} ] \label{Galatolo} Given an IET $T$ let $${\tau_r(x,y)=\min \{n\in \mathbb{N}, n>0: |T^nx-y|<r\}}.$$ For almost every IET $T$, ${\underset{ r \to 0}{\mathrm{lim\,inf}}\,\, \frac{ \log(\tau_r(x,y))}{-\log r}= 1}$ for almost every $x$.

\end{theorem}

Following result is due to L. Marchese \cite{Marchese_2011}.

\begin{theorem}[L. Marchese \cite{Marchese_2011}] \label{Marchese} Let $\{b_i\}_{i=1}^{\infty}$ be a decreasing sequence with divergent sum and with the additional property that $\{ib_i\}_{i=1}^{\infty}$ is decreasing. For almost every IET $T$ $$\delta \in \LS B(T^i(\delta'),b_i))$$ where $\delta$ and $\delta'$ are any discontinuities of $T$.
\end{theorem} Our main result is the following :

\begin{theorem}\label{Main} Let $f: [0,1) \rightarrow [0,1)$ be an ergodic measure preserving transformation (with respect to $\lambda$). Then for a decreasing sequence $\{b_i\}_{i=1}^{\infty}$ of real numbers with divergent sum

$$
\Leb(\LS f^{-i}(B (R_{\alpha}^{i} x,b_i)))=1
$$ for almost every $x \in [0,1)$ and almost every $\alpha \in [0,1)$.

\end{theorem} Corollary $\ref{Cor}$ follows directly from Theorem $\ref{Main}$.

\begin{corollary} \label{Cor} For every ergodic IET $T$ and any decreasing sequence $\{b_i\}_{i=1}^{\infty}$ of real numbers with divergent sum

$$
\lambda(\LS T^{-i}(B(R^i_{\alpha}x,b_i))) = 1
$$ for almost every $x \in [0,1)$ and almost every $\alpha \in [0,1)$.

\end{corollary}

\section{Proof of Theorem \ref{Main}}\label{proof}

Let $f: [0,1) \rightarrow [0,1)$ be an ergodic measure preserving transformation (with respect to $\lambda$). For $n\in \mathbb{N},\, m \in \mathbb{Z}$, we define sets $U_{n,m}$, $V_n$ and $W$ as follows:

\begin{equation}\label{defU}
    U_{n,m} = \{(\alpha,y) \in [0,1)\times \mathbb{R} : |f^n y -n\alpha - m| \leq b_n\} ;\,\, n\in \mathbb{N}, \, m \in \mathbb{Z}.
\end{equation}

\begin{remark}\label{rem 1} For the definition of $U_{n,m}$ to make sense we need to extend the domain of $f^n$ to $\R$. We do it as follows : Let $y \in \R$, and $\lfloor y \rfloor = y \, (\mathrm{mod} \, 1)$, then $y = \lfloor y \rfloor \, + \, k$ for some $k \in \Z$. We define $f^n (y) = f^n (\lfloor y \rfloor) + k$. With this new definition $f^n$ is well defined over entire $\mathbb{R}$.

\end{remark}
\begin{equation}
    V_n = \underset{m=-\infty}{\overset{\infty}{\bigcup}} U_{n,m}. 
\end{equation}

\begin{equation}
    W = \underset{s=1}{\overset{\infty}{\bigcap}} \underset{n=s}{\overset{\infty}{\bigcup}}V_n.
\end{equation}

\begin{proposition} Let $\{b_n\}_{n=1}^{\infty}$ be a decreasing sequence  of positive real numbers with divergent sum, then there exists a decreasing sequence $\{b'_n\}_{n=1}^{\infty}$ (of positive real numbers with divergent sum) such that $\underset{n \rightarrow \infty}{\mathrm{lim}}\, \dfrac{b'_n}{b_n} = 0$ and $b'_n < \dfrac{1}{16\,n}$ for every $n \in \mathbb{N}$. 

\end{proposition}

\noindent\textit{Proof}. We construct a sequence $\{c_n\}_{n=1}^{\infty}$ as follows : Since $\underset{n=1}{\overset{\infty}{\sum}} b_n =\infty$, there exists $n_2 \in \mathbb{N}$ such that $\underset{n=1}{\overset{n_2}{\sum}} b_n \geq 2$, define $c_n = \dfrac{b_n}{2}$ for $1 \leq n\leq n_2$. Similarly, there exists $n_3 \in \mathbb{N}$ such that $\underset{n = n_2 + 1}{\overset{n_3}{\sum}} b_n \geq 3$, define $c_n = \dfrac{b_n}{3}$ for $n_2 + 1 \leq n\leq n_3$ and so on. Put $b'_n = \mathrm{min}\,\, \{\dfrac{1}{16\,n}\,,\, c_n\}$ for $n \in \mathbb{N}$, then the sequence $\{b'_n\}_{n=1}^{\infty}$ satisfies the required conditions. \hfill{$\square$}

By replacing $\{b_n\}$ with $\{b'_n\}$ we define $U'_{n,m}$, $V'_n$ and $W'$ in the same manner as above. Let $S$ be the unit square $S = \{(\alpha,y) : 0\leq \alpha <1,0\leq y < 1\}$ and $K$ be the strip $ K = \{(\alpha,y) :0\leq \alpha < 1, y \in \mathbb{R}\}$. We denote by $\lambda_2$ the Lebesgue measure on $(0,1]^2$. 

\begin{proposition}\label{area}

For every $n \in \mathbb{N}$ 

\begin{equation*}
\lambda_2(S \cap V'_{n}) = \lambda_2 (K \cap U'_{n,0}) = 2 b'_n.
\end{equation*}

\end{proposition}

\noindent\textit{Proof}. Observe that the set $V'_{n}$ consists of countable number of copies of set $U'_{n,0}$, shifted by integer values in the $y$ direction. Hence the area $K \cap U'_{n,0}$, (where $K$ is the vertical strip in $\alpha \times y$ plane as defined above) is equal to the area $S \cap V'_{n}$, where $S$ is the unit square in the $\alpha \times y$ plane). 

For $n \in \mathbb{N}$, $m \in \mathbb{Z}$, we define following sets: 

\begin{equation}
    P'_{n,m} = \{(\alpha,f^n y) \in [0,1) \times \mathbb{R} : |f^n y -n\alpha - m| \leq b'_n\} ;\,\, n\in \mathbb{N}, \, m \in \mathbb{Z}.
\end{equation}

\begin{equation}
    Q'_n = \underset{m=-\infty}{\overset{\infty}{\bigcup}} P'_{n,m}. 
\end{equation} Denote by $z_y = f^n y$, for all $y \in \mathbb{R}$. Thus we can write
 
 \begin{equation}\label{eq defP}
P'_{n,m} = \{(\alpha,z_y) \in [0,1) \times \mathbb{R} : |z_y -n\alpha - m| < b'_n\}; \,\,\, Q'_n = \underset{m=-\infty}{\overset{\infty}{\bigcup}} P'_{n,m}.
 \end{equation} We will show

\begin{equation}
\lambda_2(S_z \cap Q'_{n}) = \lambda_2 (K_z \cap P'_{n,0}) = 2 b'_n,
\end{equation} where $S_z$ is the square $S_z = \{(\alpha,z_y) : 0\leq \alpha <1,0\leq z_y < 1\}$ and $K_z$ is the strip $K_z = \{(\alpha,z_y) :0\leq \alpha < 1, z_y \in \mathbb{R}\}$ on the $\alpha \times z$ plane. Note that $P'_{n,0} = \{(\alpha,z_y) : |z_y -n\alpha| < b'_n\}$.  Thus $\lambda_2 (K_z \cap P'_{n,0})$ is the area of the region $|z_y -n\alpha| < b'_n$, $\alpha \in [0,1)$. This area is given by 

$$
\int_{0}^{1} (b'_n + n \alpha) - (n \alpha -b'_n) \, d \alpha = 2b'_n
$$ Since $f^n$ is measure preserving $\lambda_2(S \cap V'_{n}) = \lambda_2(S_z \cap Q'_{n})$. \hfill{$\square$}

\begin{proposition}\label{parallelograms} For $j,k \in \mathbb{N}$ such that  $j<k$ we have 

\begin{equation*}
    \lambda_2 (S \cap V'_j \cap V'_k) = 4 b'_{j} b'_{k}.
\end{equation*}

\end{proposition}

\noindent\textit{Proof}. Observe that for $j < k$,

\begin{equation}\label{eq ints}
 \lambda_2(S_z \cap Q'_j \cap Q'_k) = \lambda_2 ((K_z \cap P'_{k,0})\cap Q'_j).
\end{equation} Since $f^n$ is measure preserving $(\ref{eq ints})$ is equivalent to

\begin{equation}\label{eq ints2}
 \lambda_2 (S \cap V'_j \cap V'_k) = \lambda_2((K \cap U'_{k,0})\cap V'_j). 
\end{equation} Again put $z_y = f^n y$, for all $y \in \mathbb{R}$ and consider the definition of $P'_{n,m}$ as given in $(\ref{eq defP})$. Observe that set $(K_z \cap P'_{k,0})\cap Q'_j$ consists of $k-j$ parallelograms each of area equal to $\dfrac{4b'_j b'_k}{(k-j)}$. Thus 

\begin{equation}
\lambda_2 ((K_z \cap P'_{k,0})\cap Q'_j) = 4 b'_j b'_k.
\end{equation} The equivalence of $(\ref{eq ints})$ and $(\ref{eq ints2})$ completes the proof. \hfill{$\square$} 
\begin{lemma}\label{pos. measure} Let $t \in \mathbb{N}$ be a fixed positive integer. Then for every $N_0 \in \mathbb{N}$, there exists a positive integer $N$ such that 

$$
\lambda_{2}\Big(S \cap \Big( \underset{n=N_0}{\overset{N}{\bigcup} \, V'_{tn}}\Big) \Big) > \frac{1}{8}.
$$

\end{lemma}

\noindent\textit{Proof}. Note the following inequality: 

\begin{equation}
    \lambda_{2}\Big(S \cap \Big( \underset{n=N_0}{\overset{N}{\bigcup} \, V'_{tn}}\Big) \Big) \geq \underset{n=N_0}{\overset{N}{\sum}} \lambda_2 (S \cap V'_{tn}) \,\, - \underset{N_0\leq j < k \leq N}{\sum} \lambda_2 (S \cap V'_{tj} \cap V'_{tk}).
\end{equation} By Proposition $\ref{area}$ and Proposition $\ref{parallelograms}$ for every $t \in \mathbb{N}$, 

\begin{equation*}
\lambda_2(S \cap V'_{tn}) = 2 b'_{tn} \,;\,\,\,\mathrm{and}\,\, \lambda_2 (S \cap V'_{tj} \cap V'_{tk}) = 4\,\, b'_{tj}\,\, b'_{tk}.
\end{equation*} Which implies

\begin{equation*}
    \lambda_{2}\Big(S \cap \Big( \underset{n=N_0}{\overset{N}{\bigcup} \, V'_{tn}}\Big) \Big) \geq 2 \underset{n=N_0}{\overset{N}{\sum}}b'_{tn} - \,\,4 \underset{N_0\leq j < k \leq N}{\sum} b'_{tj} \,\, b'_{tk} \geq 2 \underset{n=N_0}{\overset{N}{\sum}}b'_{tn} - \Big(2 \underset{n=N_0}{\overset{N}{\sum}} b'_{tn} \Big)^2.
\end{equation*}  Note that that for a fixed $t \in \Z$,  $\underset{n=0}{\overset{\infty}{\sum}}b'_{tn} = \infty$. To see this, assume by contradiction that there exists $c \in \mathbb{R}$, such that $\underset{n=0}{\overset{\infty}{\sum}}b'_{tn} \leq c$. Using the fact that $\{b'_n\}$ is a decreasing sequence we get, $$
\underset{n=t}{\overset{\infty}{\sum}}b'_{n} \leq t b'_{t}+ t b'_{2t} + t b'_{3t} +...+ = t \, \underset{n=1}{\overset{\infty}{\sum}}b'_{tn} \leq t\,c
$$ which is a contradiction since $\{b'_n\}$ has divergent sum. 

By definition for $n \in \mathbb{N}$, $b'_n < \dfrac{1}{16 \,n}$. Thus, $b'_{tn} < \dfrac{1}{16\,tn} \leq \dfrac{1}{16} \,\, \textrm{for}\,\, t, n \in \mathbb{N}.$ Hence we can choose the index $N$ such that

$$
\frac{3}{8} < 2 \underset{n=N_0}{\overset{N}{\sum}}b'_{tn} < \frac{3}{8} + 2.\frac{1}{16} = \frac{1}{2}. 
$$ Thus we have

$$
\lambda_{2}\Big(S \cap \Big( \underset{n=N_0}{\overset{N}{\bigcup} \, V'_{tn}}\Big) \Big) > \frac{3}{8} - (\frac{1}{2})^2 > \frac{3}{8} - \frac{1}{4} = \frac{1}{8}. 
$$ \hfill{$\square$} 

\begin{proposition}\label{Prop alpha} 

For almost all $\alpha \in [0,1)$ and a positive measure set of $x \in [0,1)$

\begin{equation}
    \Leb(\LS f^{-i}(B (R_{\alpha}^{i} x,b'_i)))> 0 . 
\end{equation}
\end{proposition} Before we prove Proposition $\ref{Prop alpha}$ we discuss the following lemma: 

\begin{lemma}\label{interval}
Let $B \subset [0,1)$ be a set such that $\lambda (B) > \sigma$ for some $0 < \sigma < 1$. Assume that there exists $k \in \mathbb{N}$ such that 

$$
B = B + \dfrac{1}{k} \,\, (\mathrm{mod} \,\,1)
$$ Let $J \subset [0,1)$ be an interval such that $\dfrac{100}{k} < l(J) < \dfrac{101}{k}$, where $l(J)$ denotes the length of $J$. Then there exists a constant $C$ (independent of $J$ and $k$) such that 

$$
\lambda(J \cap B) > C \,.\, \sigma\,.\, l(J).
$$

\end{lemma}

\noindent \textit{Proof.} Let $x \in B$ then there exists $n(x) \in \mathbb{N}$, such that $ y_x : = x + \dfrac{n(x)}{k} \,\, (\mathrm{mod} \,\,1) \in J$. This is true since $l(J)$ is much larger than $\dfrac{1}{k}$. Since $l(J) > \dfrac{100}{k}$, $y_x$ is not unique, in fact there are at least $99$ such points (taking into consideration end points of $J$). Define

$$
D_x = \{y_x \in J : y_x = x + \dfrac{n}{k},\, \textrm{for some} \,\, n \in \mathbb{N}\}.
$$ As discussed above $|D_x|$ is at least $99$. Since $B$ is invariant under $\dfrac{1}{k}$ we have $D_x \subset B$. Observe that

$$B = \underset{j=0}{\overset{k-1}{\bigsqcup}} \big(\, B \cap [0,\frac{1}{k})\,\big)+ \frac{j}{k}$$ Thus we get

$$
\lambda(B \cap J) = \int_{B \cap\, [0,\frac{1}{k})} |D_x| \,\, dx
$$

$$
\geq 99. \,\,\, \lambda \big(B \cap [0,\frac{1}{k})\big) \geq 99. \,\, \dfrac{\sigma}{k} =  \dfrac{99}{101} \,\,\sigma\,\, \dfrac{101}{k} >  \dfrac{99}{101} \,\,\sigma\,\, l(J).
$$ Put $C =  \dfrac{99}{101}$, we get

$$
\lambda(B \cap J) > C \,.\, \sigma \,.\, l(J). 
$$ \hfill{$\square$} 

\noindent\textit{Proof of Proposition \ref{Prop alpha}}. Note that by definition for $t\in \mathbb{N}$,

\begin{equation}
    U'_{tn,m} = \{(\alpha,y) \in [0,1)\times \mathbb{R} : |f^{tn} y -tn\alpha - m| \leq b'_{tn}\} ;\,\, n\in \mathbb{N}, \, m \in \mathbb{Z}.
\end{equation}

\begin{equation}
    V'_{tn} = \underset{m=-\infty}{\overset{\infty}{\bigcup}} U'_{tn,m}. 
\end{equation} For $t\in \mathbb{N}$ we define

\begin{equation}\label{eq 3a}
    W'_t = \underset{s=1}{\overset{\infty}{\bigcap}} \underset{n=s}{\overset{\infty}{\bigcup}}V'_{tn}.
\end{equation} By Lemma $\ref{pos. measure}$ we have $\lambda_2 (S \cap W'_t) > 0$. If we think of $U'_{tn,m}$ as

\begin{equation*}
    U'_{tn,m} = \{(\alpha,y) \in [0,1)\times \mathbb{R} : |f^{tn} y -tn\alpha - m| \leq b'_{tn}\} ;\,\, n\in \mathbb{N}, \, m \in \mathbb{Z}.
\end{equation*} Then $\lambda_2 (S \cap W'_t) > 0$ implies that for every $t\in \mathbb{N}$, for a positive measure set of $\alpha \in [0,1)$ and a positive measure set of $x \in [0,1)$

\begin{equation}
    \Leb(\LS f^{-ti}(B (R_{\alpha}^{ti} x,b'_{ti})))> 0 . 
\end{equation}

For $t \in \mathbb{N}$ we call such positive set of $\alpha \in [0,1)$ by $A_t$. Thus in particular for $\alpha \in A_1$, we have 

\begin{equation}\label{eq 14}
    \Leb(\LS f^{-i}(B (R_{\alpha}^{i} x,b'_{i})))> 0 . 
\end{equation} for a positive measure set of $x \in [0,1)$. Note that for every $t \in \mathbb{N}$ the set $A_t$ is invariant under addition by $\dfrac{1}{t} \,\, (\mathrm{mod} \,1)$; and $A_{t} \subseteq A_1$. Lemma \ref{pos. measure} implies that for every $t \in \mathbb{N}$, $\lambda(A_t) > \dfrac{1}{8}$. We want to show that $A_1$ has full measure.

For any $2r \in (0,1)$, we can find appropriate $t_r \in \mathbb{N}$ such that $\dfrac{100}{t_r} < 2r < \dfrac{101}{t_r}$. For any $x \in (0,1]$ by using Lemma $\ref{interval}$  we get

$$
\lambda(A_{t_r} \cap B(x, r)) > C \,.\, \dfrac{1}{8}\,.\, 2r
$$ where $B(x, r)$ denotes a ball of radius $r$ centered at $x$. Since $A_{t_r} \subset A_1$ we get

\begin{equation}\label{A_1}
    \lambda (A_1 \cap B(x, r)) > C \,.\, \dfrac{1}{8}\,.\, 2r \,=\, \dfrac{C}{4} \, r.
\end{equation} Assume $\lambda(A_1) < 1$ and let $x$ be a Lebesgue density point for $A_1^c$ then

$$\underset{r \rightarrow 0}{\mathrm{lim}} \,\, \dfrac{\lambda(A_1^c \cap B(x, r))}{\lambda(B(x, r))} = 1.$$ Thus for every $\epsilon > 0$ we can find $r > 0$ such that

$$
\lambda(A_1^c \cap B(x, r)) > (1-\epsilon)\, 2r.
$$ Choose $\epsilon = \dfrac{C}{8}$ (for $C$ as above) we get 

$$
\lambda(A_1^c \cap B(x, r)) > (1-\dfrac{C}{8})\, 2r = 2r - 2r . \dfrac{C}{8} = 2r - \dfrac{C}{4}\, r.
$$ Hence we get,

$$
\lambda(A_1 \cap B(x, r)) <  \dfrac{C}{4}\, r.
$$ which is a contradiction thus $\lambda(A_1) = 1$. This completes the proof of Proposition $\ref{Prop alpha}$.  \hfill{$\square$} 

\begin{lemma}\label{lem x&y}
For almost all $\alpha \in [0,1)$ and a.e $x \in [0,1)$

\begin{equation*}
    \Leb(\LS f^{-i}(B (R_{\alpha}^{i} x,b'_i)))= 1 . 
\end{equation*}

\end{lemma}

\noindent\textit{Proof}. So far we have shown that for almost every $\alpha \in [0,1)$, and a positive measure set of $x \in [0,1)$

\begin{equation}\label{eq pos}
    \Leb(\LS f^{-i}(B (R_{\alpha}^{i} x,b'_i)))> 0 . 
\end{equation} Fix an irrational $\alpha \in [0,1)$ from the full measure set of $\alpha$'s that satisfies $(\ref{eq pos})$ and define the set $H_{\alpha}$ as follows : $$
H_{\alpha} = \{(x,y)\in [0,1)^2 : y \in \LS f^{-i}(B (R_{\alpha}^{i} x,b'_i))\}
$$ Note that $(\ref{eq pos})$ implies that $\lambda_2(H_{\alpha}) >0 $ where $\lambda_2$ is the Lebesgue measure on $[0,1)^2$. We show that if $(x,y) \in H_{\alpha}$ then $(R_{\alpha}x,f(y)) \in H_{\alpha}$. To see this, let $(x,y) \in H_{\alpha}$ then

\begin{center}
    $f^iy$ belongs to an infinite number of sets $B (R_{\alpha}^{i} x,b'_i)$.
\end{center} Since $b'_i \leq b'_{i-1}$

\begin{center}
    $f^iy$ belongs to an infinite number of sets $B (R_{\alpha}^{i} x,b'_{i-1})$.
\end{center} Thus we get

\begin{center}
    $f^{i-1}f(y)$ belongs to an infinite number of sets $B (R_{\alpha}^{i-1} (R_{\alpha} x),b'_{i-1})$.
\end{center} Hence

\begin{center}
    $f(y)$ belongs to an infinite number of sets $f^{-(i-1)} B (R_{\alpha}^{i-1} (R_{\alpha} x),b'_{i-1})$.
\end{center} In other words

\begin{center}
    $f(y)$ belongs to an infinite number of sets $f^{-i} B (R_{\alpha}^{i} (R_{\alpha} x),b'_{i})$
\end{center} which implies $(R_{\alpha}x,f(y)) \in H_{\alpha}$.

Now we mention a condition for ergodicity of cartesian product of two ergodic measure preserving transformations of probability space. We refer our readers to \cite[Theorem 12]{Quas_2012} for a detailed proof. In the theorem below $U_T$ (resp. $U_S$) denotes the Koopman operator associated with the transformation $T$ (resp. $S$).

\vspace{2mm}

\noindent \textsc{Theorem} (Ergodicity of product). Let $T$ and $S$ be two probability measure preserving transformations. If $T$ and $S$ are ergodic, then $T \times S$ is ergodic if and only if $U_S$ and $U_T$ have no common eigenvalues other than $1$.

\vspace{2mm}

Since $f$ is a probability measure preserving transformation $U_f$ has at most countably many eigenvalues. Hence $U_f$ shares non-trivial eigenvalue with $R_{\alpha}$ for at most countably many $\alpha$. Thus by the result above $R_{\alpha} \times f$ is ergodic for almost every $\alpha \in [0,1)$. This implies $\lambda_2 (H_{\alpha}) = 1$ for almost every $\alpha \in [0,1)$. This completes the proof of the lemma.  \hfill{$\square$} 

\vspace{3mm}

\noindent\textit{Proof of Theorem $\ref{Main}$}. Observe that $b'_i < b_i$, for every $i \in \mathbb{N}$, hence

$$
\LS (B (R_{\alpha}^{i} x,b'_i)) \subset \LS (B (R_{\alpha}^{i} x,b_i)).
$$ Thus we get,

$$
\LS f^{-i}(B (R_{\alpha}^{i} x,b'_i)) \subset \LS f^{-i}(B (R_{\alpha}^{i} x,b_i)).
$$ Hence Lemma $\ref{lem x&y}$ implies Theorem $\ref{Main}$. \hfill{$\square$} 

\vspace{3mm}

\textbf{Acknowledgments.} I would like to thank Jon Chaika for suggesting this problem and many helpful discussions. I would also like to thank the referee for careful reading and valuable suggestions.

\bibliographystyle{unsrt}
\bibliography{references1.bib}

\end{document}